\documentclass{emsprocart}

\contact[aioana@ucsd.edu]{Adrian Ioana, Department of Mathematics, UC San Diego, La Jolla, CA, 92093, USA}

\newtheorem{theorem}{Theorem}[section]

\newtheorem{proposition}[theorem]{Proposition}

\theoremstyle{definition}
\newtheorem{definition}[theorem]{Definition}

\newtheorem{example}[theorem]{Example}

\title[Classification and rigidity for von Neumann algebras]{Classification  and rigidity for \\ von Neumann algebras}

\author[Adrian Ioana]{Adrian Ioana}
\begin{document}

\begin{abstract}
We survey some recent progress on the classification of von Neumann algebras arising from countable groups and their measure preserving actions on probability spaces. In particular, we present results which provide classes of  (W$^*$-superrigid) groups and actions that can be entirely reconstructed from their von Neumann algebras. We also discuss the recent finding of several large families of II$_1$ factors that have a unique group measure space decomposition.
\end{abstract}



\maketitle

\section{Introduction}

\subsection{The classification problem for II$_1$  factors}
A {\it von Neumann algebra} is an algebra of bounded linear operators on
a Hilbert space  which is closed under the $*$-adjoint operation and in the weak operator
topology.
The original papers of Murray and von Neumann \cite{MvN36, MvN43} show that countable groups
and their actions on probability spaces  provide a natural source of examples of  von Neumann algebras.
If $\Gamma$ is a countable group, then the left regular representation of $\Gamma$ on $\ell^2(\Gamma)$ generates the {\it group von Neumann algebra} $L(\Gamma)$. Every measure preserving 
action $\Gamma\curvearrowright (X,\mu)$ of a countable group $\Gamma$ on a probability space $(X,\mu)$ gives rise to the
{\it group measure space von Neumann algebra} $L^{\infty}(X)\rtimes\Gamma$. 

A central problem in the theory of von Neumann algebras is to classify $L(\Gamma)$ in the terms of the group $\Gamma$ and $L^{\infty}(X)\rtimes\Gamma$ in terms of the group action $\Gamma\curvearrowright (X,\mu)$. This problem is the most interesting in the case when $\Gamma$ has infinite
conjugacy classes (icc) and, respectively, when the action $\Gamma\curvearrowright (X,\mu)$ is free and ergodic.
These assumptions imply that the associated algebras are II$_1$ {\it factors}, i.e. 
infinite dimensional von Neumann algebras that admit a positive trace and have trivial center. Moreover, it follows that $L^{\infty}(X)$ is a {\it Cartan subalgebra} of $L^{\infty}(X)\rtimes\Gamma$, i.e. a maximal abelian von Neumann subalgebra whose normalizer generates  $L^{\infty}(X)\rtimes\Gamma$. 

The classification of group measure factors is closely related to the study of group actions up to orbit equivalence in ergodic theory. This connection is due to Singer \cite{Si55} who showed that the isomorphism class of $L^{\infty}(X)\rtimes\Gamma$ only depends on the equivalence relation given by the orbits of $\Gamma\curvearrowright (X,\mu)$. Orbit equivalence theory, initiated by Dye \cite{Dy58} and further developed by Zimmer \cite{Zi84}, has recently witnessed a flurry of activity (see the surveys \cite{Sh04,Po06b,Fu09,Ga10}). 
Over the years, the study of II$_1$ factors has also provided fertile ground for the development  of new, exciting theories:
Jones'  theory of subfactors and Voiculescu's free probability theory.

\vskip 0.1in
The goal of this survey is to present some of the spectacular progress made recently (2001-) on the classification of II$_1$ factors. 
For more on this, see the surveys \cite{Po06b,Va10a} which emphasize earlier advances.
In order to provide some context and motivation, we first give a brief history of prior work in this area. In the rest of the introduction we give a chronological account of the developments during the last decade, starting with Popa's deformation/rigidity theory. We emphasize the surge of activity since 2007 that has led to the first classes of  II$_1$ factors with a unique Cartan subalgebra (up to unitary conjugacy) and the first families of W$^*$-superrigid groups and group actions.

 \subsection{Amenable vs. non-amenable}
 Early work on the classification problem culminated with Connes' celebrated theorem \cite{Co75b} from the mid 70s:
all II$_1$ factors arising from icc amenable groups and free ergodic actions of infinite
amenable groups are isomorphic to the hyperfinite II$_1$ factor. This shows that amenable groups manifest a remarkable lack of rigidity: any algebraic property 
of the group (e.g. being torsion free) and any dynamical property of their actions (e.g. being mixing)  is lost in the passage to von Neumann algebras.

On the other hand, it was realized over the years that the non-amenable case is far more complicated and that its study underlies a beautiful rigidity theory. 
Thus, non-amenable group were used in \cite{Mc69} and \cite{Co75a}  to construct large families of non-isomorphic II$_1$ factors. In the early 80s, Connes discovered the first rigidity phenomena for von Neumann algebras  \cite{Co80}. He showed that icc property (T) groups give rise to II$_1$ factors with countable fundamental and outer automorphism groups. Representation theoretic properties of groups (property (T), Haagerup's property, weak amenability) were then used to distinguish between II$_1$ factors associated with certain lattices in Lie groups \cite{CJ83,CH88}. But, despite these rigidity results, for  a long time, the classification problem for non-amenable  II$_1$ factors remained intractable.

\subsection{Popa's deformation/rigidity theory}
A major breakthrough in the classification of II$_1$ factors was made between 2001-2004 by Popa through the invention of deformation/rigidity theory (\cite{Po01a,Po01b,Po03,Po04}, see the surveys \cite{Po06b,Va06a}). Popa discovered that if a II$_1$ factor $M$ has both a {\it deformation} and a {\it rigidity} property, then the combination of the two can sometimes be used to prove structural results for $M$ and its subalgebras.

As an illustration of this principle,  Popa proved that any II$_1$ factor $M$ admits at most one Cartan subalgebra $A$ with the relative property (T) such that  $M$ has a certain compact approximation property relative to $A$ \cite{Po01b}.  Therefore, if such a subalgebra $A$ exists, then any invariant of $M$ is equal to the corresponding invariant of the inclusion $A\subset M$. When  applied to $M=L^{\infty}(\mathbb T^2)\rtimes$ SL$_2(\mathbb Z)$ this implies
that the fundamental group of $M$ is equal to the fundamental group of the equivalence relation  of the action SL$_2(\mathbb Z)\curvearrowright\mathbb T^2$, which is trivial by Gaboriau's work \cite{Ga99,Ga01}. It altogether follows that $M$ has trivial fundamental group, thereby solving the long standing problem of finding such II$_1$ factors \cite{Po01b}. 

In \cite{Po01a,Po03}, Popa discovered that Bernoulli actions satisfy a powerful deformation property, called {\it malleability} (see Section \ref{malle}).
This allowed him to greatly expand the scope of  deformation/rigidity theory by proving the first strong rigidity theorem for group measure space factors \cite{Po03,Po04}. More precisely, let $\Gamma\curvearrowright (X,\mu)$ be a Bernoulli action of an icc group $\Gamma$ and $\Lambda\curvearrowright (Y,\nu)$ be a free ergodic probability measure preserving (pmp) action of a property (T) group $\Lambda$. It is shown in \cite{Po03,Po04} that if the group measure space factors  $L^{\infty}(X)\rtimes\Gamma$ and $L^{\infty}(Y)\rtimes\Lambda$
 are isomorphic, then the groups are $\Gamma$ and $\Lambda$ are isomorphic and their actions are {\it conjugate}. 
This means that there exist an isomorphism of probability spaces $\theta:X\rightarrow Y$ and an isomorphism of groups $\delta:\Gamma\rightarrow\Lambda$ such that $\theta(g\cdot x)=\delta(g)\cdot\theta(x)$, for all $g\in\Gamma$ and almost all $x\in X$. 

\subsection{Calculation of invariants and indecomposability for II$_1$ factors}
The classification of II$_1$ factors is very much related to the calculation of their invariants  and the study of their various decompositions. 
In the last ten years, a lot of progress has been made in both of these directions.  Recall that the {\it fundamental group} of a II$_1$ factor $M$ is the multiplicative group of $t>0$ for which the amplification $M^t$ (the ``$t$ by $t$ matrices over $M$") is isomorphic to $M$.

The calculation of the  fundamental group $\mathcal F(M)$ and the outer autormorphism group Out$(M)$ of a II$_1$ factor $M$ was and is an extremely challenging problem.  
By \cite{Co80},  these groups are countable if $M$ comes from a property (T) group. Free probability theory was later used to prove that $\mathcal F(L(\mathbb F_{\infty}))=\mathbb R_{+}^*$ \cite{Vo89,Ra91}. 

But, it was not until  Popa's deformation/rigidity theory that the first examples of II$_1$ factors with trivial fundamental and outer automorphism group were found. After giving the first examples of II$_1$ factors with trivial fundamental group \cite{Po01b}, Popa showed that $\mathcal F(M)$ can be any countable subgroup of $\mathbb R_{+}^{*}$ \cite{Po03}. 

The first calculations of Out$(M)$ were obtained in joint work with Peterson and Popa  \cite{IPP05}. We showed that for any compact abelian group $K$,  there exists a II$_1$ factor $M$ such that Out$(M)=K$. In particular, there are II$_1$ factors without outer automorphisms. The techniques of \cite{IPP05} were further exploited to show that there are II$_1$ factors with no non-trivial subfactors of finite Jones index \cite{Va06b} and that Out$(M)$ can be any compact group \cite{FV07}.

Subsequently, our understanding of the possible values of $\mathcal F(M)$ and Out$(M)$ has considerably improved.  Remarkably, Popa and Vaes showed that $\mathcal F(M)$  can be uncountable and  $\not=\mathbb R_{+}^*$, and moreover that $\mathcal F(M)$ can have arbitrary Hausdorff dimension $\alpha\in [0,1]$ (\cite{PV08a,PV08b}, see also \cite{De10}). Additionally, they proved  that Out$(M)$ can be any countable group \cite{PV06,Va07}. 
In spite of all this progress on $\mathcal F(M)$ and Out$(M)$, not a single calculation of the endomorphism semigroup End$(M)$ is yet available. Nevertheless, a description of all endomorphisms of certain II$_1$ factors $M$ was obtained in \cite{Io10}.

In the mid 90s, Voiculescu used his free entropy dimension to show that the free group factors $L(\mathbb F_n)$ do not have Cartan subalgebras  \cite{Vo95}.  Ge then showed that the free group factors are also  {\it prime} \cite{Ge96}, i.e. they cannot be written as the tensor product of two II$_1$ factors. In the last decade, these indecomposability results have been generalized and strengthened in many ways.

Using subtle C$^*$-algebraic techniques, Ozawa proved that II$_1$ factors arising from  icc hyperbolic groups $\Gamma$ are {\it solid}: the relative commutant $A'\cap L(\Gamma)$ of any diffuse von Neumann subalgebra $A\subset L(\Gamma)$ is amenable \cite{Oz03}. In particular,  $L(\Gamma)$ and all of its non-hyperfinite subfactors are prime.
Techniques  from \cite{Po01b,Oz03} were then combined in \cite{OP03} to provide a family of II$_1$ factors that can be uniquely written as a tensor product of prime factors. 
By developing a novel technique based on closable derivations, Peterson was able to show that II$_1$ factors arising from icc  groups with positive first $\ell^2$-Betti number are prime \cite{Pe06}. A new proof of solidity of $L(\mathbb F_n)$ was found by Popa in \cite{Po06c}, while II$_1$ factors coming from icc groups $\Gamma$ admitting a proper cocycle into $\ell^2(\Gamma)$ were shown to be solid in \cite{Pe06}. 
For further examples of prime and solid II$_1$ factors, see \cite{Oz04,Po06a,CI08,CH08}.

In \cite{OP07}, Ozawa and Popa discovered that the free group factors enjoy a remarkable structural property, called {\it strong solidity}, which strengthens both solidity and absence of Cartan subalgebras: the normalizer of any diffuse amenable subalgebra $A\subset L(\mathbb F_n)$ is amenable. Recently, Chifan and Sinclair showed that, more generally, the group von Neumann algebra of any icc hyperbolic group is strongly solid \cite{CS11}.
For more examples of strongly solid factors, see  \cite{OP08,Ho09,HS09,Si10}.

\subsection{W$^*$-superrigidity and uniqueness of Cartan subalgebras} Two free ergodic pmp actions $\Gamma\curvearrowright (X,\mu)$ and $\Lambda\curvearrowright (Y,\nu)$ are {\it orbit equivalent} (OE) if there exists an isomorphism of probability spaces $\theta:X\rightarrow Y$ such that $\theta(\Gamma\cdot x)=\Lambda\cdot\theta(x)$, for almost every $x\in X$.   
Singer showed that orbit equivalence amounts to the existence of an isomorphism $L^{\infty}(X)\rtimes\Gamma\cong L^{\infty}(Y)\rtimes\Lambda$ which identifies the Cartan subalgebras $L^{\infty}(X)$ and $L^{\infty}(Y)$ \cite{Si55}. Thus, W{\it $^*$-equivalence} of actions (imposing isomorphism of their group measure space factors) is weaker than orbit equivalence. In addition, orbit equivalence is clearly weaker than conjugacy.
 
By \cite{OW80}, if $\Gamma$ and $\Lambda$ are infinite amenable groups, then any free ergodic pmp actions $\Gamma\curvearrowright (X,\mu)$ and $\Lambda\curvearrowright (Y,\nu)$  are orbit equivalent. On the other hand, the culmination of a series of works \cite{GP03,Io06b,GL07,Ep07} has recently resulted in showing that any non-amenable group $\Gamma$ admits uncountably many non-OE free ergodic pmp actions \cite{Ep07} (see  the survey \cite{Ho11}).
 Furthermore, as shown in \cite{Io06b}, $\Gamma$ admits uncountably many non-W$^*$-equivalent actions.

Popa's strong rigidity theorem \cite{Po04} shows that conjugacy of the actions $\Gamma\curvearrowright (X,\mu)$ and $\Lambda\curvearrowright (Y,\nu)$   can be deduced from W$^*$-equivalence of these actions if certain conditions are imposed both on the action $\Gamma\curvearrowright (X,\mu)$ and the group 
$\Lambda$. This result made it reasonable to believe that there exist actions  $\Gamma\curvearrowright (X,\mu)$ 
which are W$^*${\it-superrigid} in the following sense: an arbitrary action $\Lambda\curvearrowright (Y,\nu)$ that is W$^*$-equivalent to $\Gamma\curvearrowright (X,\mu)$ must be conjugate to it.
Specifically, Popa asked in \cite{Po04} whether Bernoulli actions of icc property (T) groups are W$^*$-superrigid.

By \cite{FM75,Fu98}, not all Cartan subalgebras come from the group measure space construction.  To distinguish the ones that do, we call them {\it group measure space Cartan subalgebras}. 
Then an action $\Gamma\curvearrowright (X,\mu)$ is W$^*$-superrigid if and only if it is OE superrigid and $L^{\infty}(X)$ is the unique group measure space Cartan subalgebra of $L^{\infty}(X)\rtimes\Gamma$, up to unitary conjugacy.  
Each of these two properties is extremely hard to establish and their combination is even more so. 

In the last 15 years, several large families families of OE superrigid actions have been found \cite{Fu98,Po05, Po06a, Ki06, Io08, PV08c, Ki09,Fu09, PS09}, including:
\begin{enumerate}
\item
the standard action SL$_n(\mathbb Z)\curvearrowright (\mathbb T^n,\lambda^n)$, for $n\geqslant 3$ \cite{Fu98}
\item Bernoulli actions $\Gamma\curvearrowright (X,\mu)=(X_0,\mu_0)^{\Gamma}$ of many groups, including property (T) groups \cite{Po05} and products of non-amenable groups \cite{Po06a}
\item arbitrary free ergodic pmp actions of certain mapping class groups \cite{Ki06}
\item free ergodic profinite actions of property (T) groups \cite{Io08} 
\end{enumerate}

Note that the actions $\Gamma\curvearrowright (X,\mu)$ from (1) and (4) are ``virtually" OE-superrigid: if a free ergodic pmp action $\Lambda\curvearrowright (Y,\nu)$  is OE to $\Gamma\curvearrowright (X,\mu)$, then the restrictions of these actions to some finite subgroups $\Gamma_0<\Gamma$ and $\Lambda_0<\Lambda$ are conjugate.

On the other hand, the first uniqueness result for Cartan subalgebras, up to unitary conjugacy, was obtained only recently (2007) by Ozawa and Popa in their breakthrough work \cite{OP07}. They proved  that II$_1$ factors $L^{\infty}(X)\rtimes\Gamma$ associated with free ergodic profinite actions of free groups $\Gamma=\mathbb F_n$ and their direct products $\Gamma=\mathbb F_{n_1}\times\mathbb F_{n_2}\times...\times\mathbb F_{n_k}$ have a unique Cartan subalgebra. This result was extended in \cite{OP08} to groups $\Gamma$ that have the complete metric approximation property and  admit a proper cocycle into a non-amenable representation. It was then shown in \cite{Pe09} that II$_1$ factors arising from  profinite actions of groups $\Gamma$ that admit an unbounded cocycle into a mixing representation but do not have Haagerup's property,   have a unique group measure Cartan subalgebra.
 However, since none of these actions was known to be OE superrigid, W$^*$-superrigidity could not be concluded. 
 
The situation changed starting with the work of Peterson \cite{Pe09} who was able to show the existence of virtually W$^*$-superrigid actions.

Shortly after, Popa and Vaes discovered the first concrete families of W$^*$-superrigid actions \cite{PV09}. They first showed a general unique decomposition result:  any II$_1$ factor arising from an arbitrary free ergodic pmp action $\Gamma\curvearrowright (X,\mu)$ of a group $\Gamma$ in a large class  of amalgamated free product groups has a unique group measure space Cartan subalgebra, up to unitary conjugacy.  Applying OE superrigidity results from \cite{Po05,Po06a} and \cite{Ki09} then allowed them to provide several wide classes of W$^*$-superrigid actions.
 For related results, see \cite{FV10,HPV10}.
 
Despite this progress, the original question of whether Bernoulli actions of icc property (T) groups are W$^*$-superrigid remained open for some time, until it was answered in the affirmative in \cite{Io10}.  The starting point of the proof is the observation that every group measure space decomposition $M=L^{\infty}(Y)\rtimes\Lambda$ of a II$_1$ factor $M$ gives rise to an embedding $\Delta:M\rightarrow M\bar{\otimes}M$ \cite{PV09}. In the case when $M$ arises from a Bernoulli action $\Gamma\curvearrowright (X,\mu)=(X_0,\mu_0)^{\Gamma}$ of an icc property (T) group $\Gamma$,  a classification of all such embeddings was given in \cite{Io10}.  
This classification is precise enough to imply that the Cartan subalgebras $L^{\infty}(X)$ and $L^{\infty}(Y)$ are unitarily conjugated.
In combination with the OE superrigidity theorem \cite{Po05} it follows that the action $\Gamma\curvearrowright (X,\mu)$ is W$^*$-superrigid.

An icc group $\Gamma$ is called W$^*${\it-superrigid} if any group $\Lambda$ satisfying $L(\Gamma)\cong L(\Lambda)$ must be isomorphic to $\Gamma$. The superrigidity question for groups is significantly harder than for actions. While large families of W$^*$-superrigid actions were found in \cite{PV09,Io10},  not a single example of a W$^*$-superrigid group was known until our joint work with Popa and Vaes \cite{IPV10}. 
We proved that  many generalized wreath product groups are W$^*$-superrigid, although plain wreath product groups essentially never have this property. For instance, we showed that given any non-amenable group $\Gamma_0$, its canonical ``augmentation" $\Gamma = (\mathbb Z/2\mathbb Z)^{(I)}\rtimes(\Gamma_0\wr\mathbb Z)$ is W$^*$-superrigid, where the set $I$ is the quotient $(\Gamma_0\wr\mathbb Z)/\mathbb Z$ on which the wreath product group $\Gamma_0\wr\mathbb Z =\Gamma_0^{(\mathbb Z)}\rtimes\mathbb Z$ acts by left multiplication.

A general conjecture predicts that II$_1$ factors $L^{\infty}(X)\rtimes\Gamma$ arising from arbitrary free ergodic pmp actions of groups $\Gamma$ with positive first $\ell^2$-Betti number, have a unique Cartan subalgebra, up to unitary conjugacy. The condition $\beta_1^{(2)}(\Gamma)>0$ is equivalent to $\Gamma$ being non-amenable and having an unbounded cocycle into its left regular representation \cite{BV97,PT07}, and is satisfied by any free product group $\Gamma=\Gamma_1*\Gamma_2$, with $|\Gamma_1|\geqslant 2$ and $|\Gamma_2|\geqslant 3$.

Several recent results provide supporting evidence for this conjecture.  Popa and Vaes showed in \cite{PV09} that if $\Gamma$ is the free product of a non-trivial group and an infinite property (T) group, then any II$_1$ factor $L^{\infty}(X)\rtimes\Gamma$ associated with a free ergodic pmp action of $\Gamma$ has a unique group measure space Cartan subalgebra. Chifan and Peterson then proved that the same holds for groups $\Gamma$ with positive first $\ell^2$-Betti number that admit a non-amenable subgroup with the relative property (T) \cite{CP10}. A common generalization of the last two results was obtained in \cite{Va10b}.  Most recently, it was proven in \cite{Io11a,Io11b} that if a group $\Gamma$ satisfies $\beta_1^{(2)}(\Gamma)>0$, then $L^{\infty}(X)\rtimes\Gamma$ has a unique group measure space Cartan subalgebra whenever the action $\Gamma\curvearrowright (X,\mu)$ is either rigid or profinite.

Generalizing \cite{OP07}, it was shown in \cite{CS11,CSU11} that profinite actions of hyperbolic groups and of direct products of hyperbolic groups give rise to II$_1$ factors with a unique Cartan subalgebra. The proofs of \cite{OP07,CS11,CSU11} rely both on the fact that free groups (and, more generally, 
hyperbolic groups, see \cite{Oz07, Oz10}) are weakly amenable and that the actions are profinite. 

In a very recent breakthrough, Popa and Vaes removed the assumption that the action is profinite and obtained wide-ranging unique Cartan subalgebra results. They proved that
if $\Gamma$ is either a weakly amenable group with $\beta_1^{(2)}(\Gamma)>0$ \cite{PV11} or a hyperbolic group \cite{PV12} (or 
a direct product of groups in one of these classes), then II$_1$ factors $L^{\infty}(X)\rtimes\Gamma$ arising from arbitrary free ergodic pmp actions of $\Gamma$ have a unique Cartan subalgebra, up to unitary conjugacy. In particular, this settles the ``$\beta_1^{(2)}(\Gamma)>0$ conjecture"  in the key case when $\Gamma$ is a free group. Note that in the meantime, the main result of \cite{PV12} has been extended by Houdayer and Vaes to cover non-amenable non-singular actions of hyperbolic groups  \cite{HV12}.

\vskip 0.1in
\noindent
{\bf Organization of the paper}. Besides the introduction and a section of preliminaries, this paper has four other sections. In Section 3 we review the basic notions and techniques of Popa's deformation/rigidity theory. In Sections 4-6 we expand upon the topics discussed in the last part of the introduction. In doing so, we follow three directions: W$^*$-superrigidity for Bernoulli actions and wreath product groups, uniqueness of group measure space Cartan subalgebras, and uniqueness of arbitrary Cartan subalgebras. 

\vskip 0.1in
\noindent
{\bf Acknowledgments}. It is my pleasure to thank Jesse Peterson and Sorin Popa for useful comments, and Stefaan Vaes for many suggestions that helped improve the exposition.

\section{Preliminaries} 

\subsection{Tracial von Neumann algebras}

\begin{definition} A {\it tracial von Neumann algebra} $(M,\tau)$ is a von Neumann algebra $M$ together with a positive linear functional $\tau:M\rightarrow\mathbb C$ that is
\begin{itemize}
 \item a {\it trace:} $\tau(xy)=\tau(yx)$, for all $x,y\in M$, 
 \item {\it faithful:} if $\tau(x)=0$, for some $x\geqslant 0$, then $x=0$, and
 \item {\it normal}: $\tau(\sum_{i\in I}p_i)=\sum_{i\in I}\tau(p_i)$, for any family $\{p_i\}_{i\in I}\subset M$ of mutually orthogonal projections.
 \end{itemize}

\end{definition}
Any tracial von Neumann algebra $(M,\tau)$ admits a canonical (or standard) representation on a Hilbert space. Indeed, 
denote by $L^2(M)$ the Hilbert space obtained by completing $M$ w.r.t. the 2-norm: $\|x\|_2=\sqrt{\tau(x^*x)}$. Then the left multiplication on $M$ extends to an injective $*$-homomorphism $\pi:M\rightarrow\mathbb B(L^2(M))$. 
 
 \vskip 0.1in
 \noindent
Next, we recall from \cite{MvN36,MvN43,FM75} three general constructions of tracial von Neumann algebras: from groups, group actions and equivalence relations. 

 Let $\Gamma$ be a countable group and denote by $\{\delta_h\}_{h\in\Gamma}$ the usual  orthonormal basis of  $\ell^2(\Gamma)$. The left regular representation $u:\Gamma\rightarrow\mathcal U(\ell^2(\Gamma))$ is given by $u_g(\delta_h)=\delta_{gh}$. The {\it group von Neumann algebra} $L(\Gamma)$ is the von Neumann algebra generated by $\{u_{g}\}_{g\in\Gamma}$, i.e. the  weak operator  closure of the group algebra $\mathbb C\Gamma=$ span $\{u_g\}_{g\in\Gamma}$. It has a faithful normal trace given by
$\tau(T)=\langle T\delta_e,\delta_e\rangle$. In other words, $$\tau(u_g)=\begin{cases} 1,\;\;\text{if}\;\;g=e\\ 0,\;\;\text{if}\;\; g\not=e\end{cases}$$

Now, let $\Gamma\curvearrowright (X,\mu)$ be a probability measure preserving (pmp) action of a countable group $\Gamma$ an a probability space $(X,\mu)$. Denote by $(\sigma_g)_{g\in\Gamma}$ the associated action of $\Gamma$ on $L^{\infty}(X)$, i.e. $\sigma_g(a)(x)=a(g^{-1}\cdot x).$
Then both $\Gamma$ and $L^{\infty}(X)$ are represented on the Hilbert space $L^2(X,\mu)\bar{\otimes}\ell^2(\Gamma)$ through the formulae: 

$$u_g(b\otimes \delta_h)=\sigma_g(b)\otimes\delta_{gh},\;\;\;\text{and}\;\;\; a(b\otimes \delta_h)=ab\otimes \delta_h.$$

The {\it group measure space von Neumann algebra} $L^{\infty}(X)\rtimes\Gamma$ is the von Neumann algebra generated by $\{u_g\}_{g\in\Gamma}$ and $L^{\infty}(X)$. It has a faithful normal trace defined by $\tau(T)=\langle T(1\otimes\delta_e),1\otimes\delta_e\rangle$. 
Since $u_gau_g^*=\sigma_g(a)$, the group measure space algebra $L^{\infty}(X)\rtimes\Gamma$ is equal to the weak operator closure of the span of the set $\{au_g|a\in L^{\infty}(X),g\in\Gamma\}$.
The restriction of $\tau$ to this set is given by
$$\tau(au_g)=\begin{cases} \int_{X}a\;\text{d}\mu,\;\;\text{if}\;\;g=e\\ 0,\;\;\text{if}\;\; g\not=e\end{cases}$$

Finally, let $R$ be a countable pmp equivalence relation on a probability space $(X,\mu)$ \cite{FM75}. This means that $R$ has countable classes,  $R$ is a Borel subset of $X\times X$,  and that every Borel automorphism $\theta$ of $X$ satisfying $(\theta(x),x)\in R$, almost everywhere,  preserves the measure $\mu$. The group of such automorphisms $\theta$ is called the {\it full group} of $R$ and denoted $[R]$.  Further, we consider on  $R$ the infinite Borel measure given by $$\tilde\mu(A)=\int_{X}|\{y\in X|(x,y)\in A\}|\;\text{d}\mu(x),\;\;\text{for every Borel subset}\;A\subset R.$$

Then both $[R]$ and $L^{\infty}(X)$ are represented on the Hilbert space $L^2(R,\tilde\mu)$ through the formulae $u_{\theta}(f)(x,y)=f(\theta^{-1}(x),y)$ and $(af)(x,y)=a(x)f(x,y),$
for all $f\in L^2(R,\tilde\mu)$ and every $\theta\in[R], a\in L^{\infty}(X)$. The {\it generalized group measure space von Neumann algebra} $L(R)$ associated with $R$ is generated by $\{u_{\theta}\}_{\theta\in [R]}$ and $L^{\infty}(X)$. This algebra also has a faithful normal trace given by $\tau(T)=\langle T(1_{\Delta}),1_{\Delta}\rangle$, where $\Delta=\{(x,x)|x\in X\}$, or explicitly by $$\tau(u_{\theta})=\mu(\{x\in X|\theta(x)=x\}),\;\;\;\text{and}\;\;\;\tau(a)=\int_{X}a\;\text{d}\mu(x).$$

More generally, to every scalar 2-cocycle $w$ on $R$,  one can associate a tracial von Neumann algebra $L(R,w)$, which, as above, contains a copy of $L^{\infty}(X)$ \cite{FM75}.

Note that if $\Gamma\curvearrowright (X,\mu)$ is a free pmp action and we denote by $R_{\Gamma\curvearrowright X}$ the equivalence relation given by its orbits, then there exists a canonical isomorphism  $L^{\infty}(X)\rtimes\Gamma=L(R_{\Gamma\curvearrowright X})$ whose restriction to $L^{\infty}(X)$ is the identity.

\subsection{II$_1$ factors and Cartan subalgebras}
A II$_1$ {\it factor} is an infinite dimensional tracial von Neumann algebra with trivial center.
We have that $L(\Gamma)$ is a II$_1$ factor if and only if $\Gamma$ is icc and that
 $L^{\infty}(X)\rtimes\Gamma$ is a II$_1$ factor whenever the action $\Gamma\curvearrowright (X,\mu)$ is   free and ergodic. Moreover,
 $L(R)$ is a II$_1$ factor if and only if $R$ is ergodic: any $R$-invariant Borel subset of $X$ has measure $0$ or $1$.

Let $A$ be a von Neumann subalgebra of a II$_1$ factor $M$. The {\it normalizer of $A$ in $M$}, denoted $\mathcal N_{M}(A)$, is the group of unitaries $u\in M$ satisfying $uAu^*=A$. We say that $A$ is a {\it Cartan subalgebra} of $M$ if it is maximal abelian  and its normalizer generates $M$. 
If $\Gamma\curvearrowright (X,\mu)$ is a free ergodic pmp action, then $L^{\infty}(X)$ is a {\it group measure space} Cartan subalgebra of $L^{\infty}(X)\rtimes\Gamma$. If $R$ is an ergodic pmp equivalence relation, then $L^{\infty}(X)$ is a Cartan subalgebra of $L(R)$. More generally,   $L^{\infty}(X)$ is a Cartan subalgebra of $L(R,w)$, for every 2-cocycle $w$ on $R$. Conversely, Feldman and Moore \cite{FM75} showed that any Cartan subalgebra inclusion arises in this way.

By \cite{CFW81}, any two Cartan subalgebras of the hyperfinite II$_1$ factor are conjugated by an automorphism. The first examples of II$_1$ factors admitting two Cartan subalgebras that are not conjugated by an automorphism were given by Connes and Jones in \cite{CJ81}. Examples of such II$_1$ factors where the two Cartan subalgebras are explicit were recently found by Ozawa and Popa (\cite{OP08}, see also \cite{PV09}). Very recently, a class of II$_1$ factors $M$ whose Cartan subalgebras cannot be classified, in the sense that the equivalence relation of being conjugated by an automorphism of $M$ is not Borel, has been constructed in \cite{SV11}.

Proving uniqueness of Cartan subalgebras plays a crucial role in the classification of group measure space II$_1$ factors.
As the next result shows, it allows one to reduce the classification of group measure space factors, up to isomorphism, to the classification of the corresponding group actions, up to orbit equivalence.
  
\begin{proposition} [Singer, \cite{Si55}]\label{Si} 
Let $\Gamma\curvearrowright (X,\mu)$ and $\Lambda\curvearrowright (Y,\nu)$ be free ergodic pmp actions. Then the following are equivalent:
\begin{itemize}
\item there exists a $*$-isomorphism $\theta:L^{\infty}(X)\rtimes\Gamma\rightarrow L^{\infty}(Y)\rtimes\Lambda$  such that $\theta(L^{\infty}(X))=L^{\infty}(Y)$.
\item the actions $\Gamma\curvearrowright (X,\mu)$ and $\Lambda\curvearrowright (Y,\nu)$ are orbit equivalent: there exists an isomorphism  $\alpha:X\rightarrow Y$ of probability spaces  such that $\alpha(\Gamma\cdot x)=\Lambda\cdot\alpha(x)$, for almost all $x\in X$.
\end{itemize}
\end{proposition}

The proof of this proposition relies on the fact that any unitary element $u$ in $L^{\infty}(X)\rtimes\Gamma=L(R_{\Gamma\curvearrowright X})$ that normalizes $L^{\infty}(X)$ is the form $u=au_{\theta}$, for some unitary $a\in L^{\infty}(X)$ and some $\theta\in [R_{\Gamma\curvearrowright X}]$.

\subsection{Popa's intertwining-by-bimodules}  

Let $\Gamma\curvearrowright (X,\mu)$ and $\Lambda\curvearrowright (Y,\nu)$ be free ergodic pmp actions and $\theta:L^{\infty}(X)\rtimes\Gamma\rightarrow L^{\infty}(Y)\rtimes\Lambda$ be a $*$-isomorphism. Then both $\theta(L^{\infty}(X))$ and $L^{\infty}(Y)$ are Cartan subalgebras of $M=L^{\infty}(Y)\rtimes\Lambda$. By Proposition \ref{Si}, in order to conclude that the initial actions are orbit equivalent, it suffices to find an automorphism $\rho$ of $M$ such that $\rho(\theta(L^{\infty}(X))=L^{\infty}(Y)$. 

In \cite{Po01b,Po03}, Popa developed a powerful technique for showing unitary conjugacy of subalgebras of a tracial von Neumann algebra. 
In particular, this provides a criterion for the existence of an inner automorphism $\rho$ of $M$ with the desired property.

\begin{theorem}[Popa, \cite{Po01b,Po03}]\label{inter} Let $P,Q$ be von Neumann subalgebras of a tracial von Neumann algebra $(M,\tau)$. Then the following conditions are equivalent:

\begin{itemize}
\item there exist projections $p\in P$, $q\in Q$, a $*$-homomorphism $\theta:pPp\rightarrow qQq$ and a non-zero partial isometry $v\in qMp$ such that $\theta(x)v=vx$, for all $x\in pPp$.
\item there is no sequence of unitaries $u_n\in P$ satisfying $\|E_Q(xu_ny)\|_2\rightarrow 0$, for all $x,y\in M$.
\end{itemize} 

If one of these conditions is satisfied, we say that a corner of $P$ embeds into $Q$.

Moreover, if $P$ and $Q$ are Cartan subalgebras of $M$, and a corner of $P$ embeds into $Q$, then there exists a unitary $u\in M$ such that $uPu^*=Q$.
\end{theorem}

Here,  as usual, $E_Q:M\rightarrow Q$ denotes the unique $\tau$-preserving  conditional expectation onto $Q$.

\section{Popa's deformation/rigidity theory}

\subsection{Deformations}\label{malle} In the last decade, Popa's deformation/rigidity theory has completely reshaped the landscape of von Neumann algebras.
 At the heart of Popa's   theory is the notion of deformation of  II$_1$ factors. In the first part of this section, we define this notion and then illustrate it with many examples.
 
 \begin{definition}
A {\it deformation} of the identity  of a tracial von Neumann algebra $(M,\tau)$  is a sequence $\phi_n:M\rightarrow M$ of unital, trace preserving, completely positive maps $\phi_n:M\rightarrow M$ satisfying $\|\phi_n(x)-x\|_2\rightarrow 0$, for all $x\in M$. 
\end{definition}

A linear map $\phi:M\rightarrow M$ is {\it completely positive} if for all $m\geqslant 1$ the amplification $\phi^{(m)}:\mathbb M_m(M)\rightarrow \mathbb M_m(M)$
given by $\phi^{(m)}([x_{i,j}])=[\phi(x_{i,j})]$ is positive. 

\begin{example}
Let $\phi_n:\Gamma\rightarrow\mathbb C$ be a sequence of positive definite functions on a countable group $\Gamma$ such that $\phi_n(e)=1$, for all $n$, and $\phi_n(g)\rightarrow 1$, for all $g\in\Gamma$. Then we have a deformation of $L(\Gamma)$ given by $u_g\mapsto\phi_n(g)u_g$, and one of any group measure space algebra $L^{\infty}(X)\rtimes\Gamma$ given by $au_g\mapsto \phi_n(g)au_g$.

If $\Gamma$ has Haagerup's property \cite{Ha79}, then there is such a sequence $\phi_n:\Gamma\rightarrow\mathbb C$ satisfying $\phi_n\in c_0(\Gamma)$,  for all $n$. In this case, the resulting deformation of $L^{\infty}(X)\rtimes\Gamma$ is compact relative to $L^{\infty}(X)$. This fact is a crucial ingredient of Popa's proof that the II$_1$ factor $L^{\infty}(\mathbb T^2)\rtimes\text{SL}_2(\mathbb Z)$ has trivial fundamental group \cite{Po01b}.
\end{example}

\begin{example}\label{mall} Let $\Gamma$ be a countable group and $\Gamma\curvearrowright (X,\mu)=([0,1], \text{Leb})^{\Gamma}$ be the Bernoulli action: $g\cdot (x_h)_{h\in\Gamma}=(x_{g^{-1}h})_{h\in\Gamma}$. Popa discovered in \cite{Po01a,Po03} that Bernoulli actions have a remarkable deformation property, called {\it malleability}: there exists a continuous  family of automorphisms $(\alpha_t)_{t\in\mathbb R}$ of the probability space $X\times X$, which commute with diagonal action of $\Gamma$, and satisfy $\alpha_0=\text{id}$, $\alpha_1(x,y)=(y,x).$ 

To see this, first construct a continuous family of automorphisms $(\alpha_t^0)$ of the probability space $[0,1]\times [0,1]$ such that $\alpha_0^0=\text{id}$ and $\alpha_1^0(x,y)=(y,x)$. For instance, one can take $$\alpha_t^0(x,y)=
\begin{cases} (x,y),\;\;\;\text{if}\;\;\; |x-y|\geqslant t\\ (y,x),\;\;\;\text{if}\;\;\; |x-y|<t\end{cases}$$Then identify $X\times X=([0,1]\times[0,1])^{\Gamma}$ and define   $\alpha_t((x_h)_{h\in\Gamma})=(\alpha_t^0(x_h))_{h\in\Gamma}$.

Next, let us explain how to get a deformation of $M=L^{\infty}(X)\rtimes\Gamma$ from $(\alpha_t)_{t\in\mathbb R}$. Denote by $\theta_t$ the automorphism of $L^{\infty}(X\times X)$ given by $\theta_t(a)(x)=a(\alpha_t^{-1}(x))$. Since $\theta_t$ commutes with the diagonal action of $\Gamma$,  it extends to an automorphism of $\tilde M=L^{\infty}(X\times X)\rtimes\Gamma$ by letting $\theta_t(u_g)=u_g$, for  $g\in\Gamma$. Then $(\theta_t)_{t\in\mathbb R}$ is a continuous family of automorphisms of $\tilde M$ such that $\theta_0=\text{id}$.  

Now, in general, whenever one has  such a pair $(\tilde M,(\theta_t)_{t\in\mathbb R})$, letting $$\phi_t=E_M\circ{\theta_t}_{|M}:M\rightarrow M$$ and choosing any sequence $t_n\rightarrow 0$, gives a deformation $(\phi_{t_n})_n$  of $M$. This is why from now on, abusing notation, we say that the pair $(\tilde M,(\theta_t)_{t\in\mathbb R})$ is a deformation of $M$.
\vskip 0.1in
Since they were first introduced in \cite{Po01a,Po03}, malleable deformations have been found in several other contexts and are now a central tool in  Popa's deformation/rigidity theory (see \cite[Section 6]{Po06b}) . Next, we review several constructions of malleable deformations.
 \end{example}

\begin{example}\label{wmall}
A malleable deformation for Bernoulli actions related to the one in \cite{Po03} was introduced in \cite{Io06a}. Note that if $(B,\tau)$ is a tracial von Neumann algebra and $I$ is a set, then one can construct a tensor product von Neumann algebra $\bar{\otimes}_{i\in I}B$. This algebra is  tracial, with its trace given by  $\tilde\tau(x)=\prod_{i\in I}\tau(x_i)$, for every element $x=\otimes_{i\in I} x_i$ whose support, $\{i\in I|x_i\not=1\}$, is finite.

 With the  notations from Example \ref{mall}, we can now identify $A=L^{\infty}(X)$ with the infinite tensor product algebra $\bar{\otimes}_{g\in\Gamma}A_0$, where $A_0=L^{\infty}([0,1])$.
 Define $\tilde A_0$ to be the free product von Neumann algebra $A_0*L(\mathbb Z)$.  Let  $u\in L(\mathbb Z)$ be the canonical generating unitary and choose a self-adjoint operator $h$ such that   $u=\exp(ih)$. For $t\in\mathbb R$,  let $\theta_t^0$ be the inner automorphism of $\tilde A_0$ given by $\theta_t^0=\text{Ad}(\exp(ith))$. 
Then $\theta_t=\otimes_{g\in\Gamma}\theta_t^0$ is an automorphism of $\tilde A=\bar{\otimes}_{g\in\Gamma}\tilde A_0$ which commutes with the Bernoulli action of $\Gamma$. Thus, $\theta_t$ extends to an automorphism of $\tilde M=\tilde A\rtimes\Gamma$ by letting $\theta_t(u_g)=u_g$.
\end{example}

\begin{example}\label{ipp} Let $(M_1,\tau_1)$ and $(M_2,\tau_2)$ be tracial von Neumann algebras with a common von Neumann subalgebra $A$ such that ${\tau_1}_{|A}={\tau_2}_{|A}$.  Denote by $M$ the amalgamated free product von Neumann algebra $M_1*_{A}M_2$ (see \cite{Po93} and \cite{VDN92} for the definition).  
Following \cite{IPP05}, $M$ admits a natural malleable deformation. 
More precisely, define $\tilde M=M*_A(A\bar{\otimes}L(\mathbb F_2))$. Then $M\subset\tilde M$ and one constructs a 1-parameter  group of automorphisms $(\theta_t)_{t\in\mathbb R}$ of $\tilde M$ as follows. Let $u_1$ and $u_2$ be the canonical generating unitaries of $L(\mathbb F_2)$. Choose self-adjoint operators $h_1$ and $h_2$  such that $u_1=\exp(ih_1)$ and $u_2=\exp(ih_2)$. Then  $\theta_t$ is defined to be the identity on $L(\mathbb F_2)$ and inner on both $M_1$ and $M_2$: $${\theta_t}_{|M_1}=\text{Ad}(\text{exp}(ith_1))_{|M_1} \;\;\;\text{and}\;\;\;\; {\theta_t}_{|M_2}=\text{Ad}(\text{exp}(ith_2))_{|M_2}.$$

See also \cite{Po86,Po06c}  for the construction of a malleable deformation for the free group factors. 
\end{example}

\begin{example}\label{si} Next, we recall the construction of a malleable deformation from cocycles (\cite{Si10}, see also \cite{PS09}). Start with a countable group $\Gamma$, an orthogonal representation $\pi:\Gamma\rightarrow\mathcal O(H)$ onto a separable real Hilbert space $H$, and a {\it cocycle} $c:\Gamma\rightarrow H$. In other words, the  cocycle relation $c(gh)=c(g)+\pi(g)c(h)$ holds for all $g,h\in \Gamma$. Let $\Gamma\curvearrowright (X,\mu)$ be a pmp action. Define $A=L^{\infty}(X)$ and $M=A\rtimes\Gamma$.

Let $(D,\tau)$ be the unique tracial von Neumann algebra generated by unitaries $\omega(\xi)$, with $\xi\in H$, subject to the relations $\omega(\xi+\eta)=\omega(\xi)\omega(\eta)$, $\omega(\xi)^*=\omega(-\xi)$ and $\tau(\omega(\xi))=\exp{(-\|\xi\|^2)}$. Consider the Gaussian action of $\Gamma$ on $D$ which on the generating unitaries  $\omega(\xi)$ is given by $\sigma_g(\omega(\xi))=\omega(\pi(g)(\xi))$. Further, consider the diagonal product action of $\Gamma$ on $A\bar{\otimes}D$ and denote $\tilde M=(A\bar{\otimes}D)\rtimes\Gamma$. 

Then $M\subset\tilde M$. Defining $\theta_t$ to be the identify on $A\bar{\otimes}D$ and $$\theta_t(u_g)=(1\otimes\omega(tc(g)))\;u_g\;\;\;\text{for all}\;\;\;g\in\Gamma$$
gives a 1-parameter group of automorphisms $(\theta_t)_{t\in\mathbb R}$ of $\tilde M$. 
\end{example}

\begin{example}\label{chsi} Now, assume that $c:\Gamma\rightarrow H$ is a {\it quasi-cocycle} rather than a cocycle. Thus, the cocycle relation only holds up  to bounded error: there exists  $\kappa>0$ such that $\|c(gh)-c(g)-\pi(g)c(h)\|\leqslant\kappa,$ for all $g,h\in\Gamma$. 

It was recently discovered in \cite{CS11} that quasi-cocycles can still be used to construct deformations. 
 It is clear that defining $\theta_t$ as in the cocycle case will not not work. 
The original idea of \cite{CS11} is to use the canonical unitary implementation of $\theta_t$ instead. More precisely, identify $L^2(\tilde M)$ with $L^2(A)\otimes L^2(D)\otimes\ell^2(\Gamma)$ and define a unitary operator $V_t$ on $L^2(\tilde M)$ by letting $$V_t(\xi\otimes\eta\otimes\delta_g)=\xi\otimes \omega(tc(g))\eta\otimes\delta_g.$$

For $x\in\tilde M$, define $\theta_t(x)=V_txV_t^*$. If $c$ is a cocycle, then this formula coincides with the above one and thus $\theta_t$ gives an automorphism of $\tilde M$. In general, it is only true that $\theta_t$ leaves invariant the larger von Neumann algebra $\tilde{\mathcal M}$ generated by $\tilde M$ and $\ell^{\infty}(\Gamma)$. Since $\tilde{\mathcal M}$ is not tracial, we cannot derive a deformation of $M$ in the usual sense. Nevertheless, $\theta_t$ is a ``deformation at the C$^*$-algebraic level", in a sense made precise in \cite{CS11}.
\end{example}

\begin{example} Deformations arise naturally from closable derivations  \cite{CS03,Pe06}. Let $(M,\tau)$ be a tracial von Neumann algebra and $\delta:M\rightarrow H$ be a closable, densely defined, real derivation into an $M$-$M$ bimodule $H$. Defining $\phi_t=e^{-t\delta^*\bar{\delta}}:M\rightarrow M$ we obtain a continuous semigroup $(\phi_t)_{t\geqslant 0}$ of unital trace preserving completely positive maps. It what recently shown \cite{Da10} that any such semigroup admits a dilation: there exists a tracial von Neumann algebra $\tilde M$ containing $M$ and a continuous family $(\theta_t)_{t\geqslant 0}$ of automorphisms of $\tilde M$ such that $\phi_t=E_{M}\circ{\theta_t}_{|M}$.
\end{example}

\vskip 0.2in
\subsection{Rigidity} A second central notion in  deformation/rigidity theory is Popa's relative property (T) for inclusions of von Neumann algebras.

\begin{definition} A von Neumann subalgebra $P$ of a tracial von Neumann algebra $(M,\tau)$ has the {\it relative property (T)} if any deformation $\phi_n:M\rightarrow M$ of the identity of $M$ must converge uniformly to the identity on the unital ball of $P$ \cite{Po01b}.
\end{definition}

Note that if $P=M$, then this property amounts to the property (T) of $M$, in the sense Connes and Jones \cite{CJ83}. 

Given two countable groups $\Gamma_0<\Gamma$, the inclusion of group von Neumann algebras $L(\Gamma_0)\subset L(\Gamma)$ has the relative property (T) if and only if the inclusion $\Gamma_0<\Gamma$ has the relative property (T)  \cite{Po01b}. Examples of inclusions of groups with the relative property (T) include $\text{SL}_n(\mathbb Z)<\text{SL}_n(\mathbb Z)$, for $n\geqslant 3$ \cite{Ka67}, and $\mathbb Z^2<\mathbb Z^2\rtimes\Gamma$, for any non-amenable subgroup $\Gamma<\text{SL}_2(\mathbb Z)$ \cite{Bu91}. 

Several classes of inclusions of von Neumann algebras inclusions satisfying the relative property (T) that do not arise from inclusions of  groups have been recently found in \cite{Io09,CI09,IS10}. For instance, the main result of \cite{Io09} asserts that if $M\subset L(\mathbb Z^2\rtimes\text{SL}_2(\mathbb Z))$ is a non-hyperfinite subfactor which contains $L(\mathbb Z^2)$, then the inclusion $L(\mathbb Z^2)\subset M$ has the relative property (T). \vskip 0.05in

If $(M,\tau)$ is a tracial von Neumann algebra, then an $M$-$M$ bimodule is a Hilbert space endowed with commuting $*$-representations of $M$ and its opposite algebra, $M^{\text{op}}$. Examples of $M$-$M$ bimodules include the trivial bimodule $L^2(M)$ and the coarse bimodule $L^2(M)\bar{\otimes} L^2(M)$. These  are particular instances of the bimodule $L^2(M)\bar{\otimes}_{P}L^2(M)$ obtained by  completing the algebraic tensor product $M\otimes M$ with respect to the scalar product $\langle x_1\otimes x_2,y_1\otimes y_2\rangle=\tau(y_2^*E_P(x_2^*x_1)y_1).$  Note that this bimodule is isomorphic to the $L^2$-space of Jones's basic construction $\langle M,e_P\rangle$.

The natural correspondence between completely positive maps and bimodules allows to reformulate relative property (T) in terms of bimodules. 
Recall that if $H$ is an $M$-$M$ bimodule, then a vector $\xi\in H$ is {\it tracial}  if $\langle x\xi,\xi\rangle=\langle\xi x,\xi\rangle=\tau(x)$, for all $x\in M$; it is $P$-{\it central} if $x\xi=\xi x,$ for all $x\in P$. Also, a net of vectors $(\xi_n)_n$ is called {\it almost central}  if $\|x\xi_n-\xi_n x\|\rightarrow 0$, for all $x\in M$. 
Then the relative property (T) for an inclusion $P\subset M$ means that
any $M$-$M$ bimodule without $P$-central vectors does not admit a net $(\xi_n)_n$ of tracial, almost central vectors \cite{Po01b}. 

Thus, relative property (T) requires that all bimodules satisfy a certain spectral gap property.
One way to obtain less restrictive notions of rigidity is to impose that only certain bimodules have spectral gap (see \cite[Section 6.5]{Po06b}).

 This brings us to the following:

\begin{definition} 
A tracial von Neumann algebra $(M,\tau)$ is {\it amenable} if there exists a net $\xi_n\in L^2(M)\bar{\otimes} L^2(M)$ of tracial, almost central vectors \cite{Po86}.

If $P,Q\subset M$ are von Neumann subalgebras,  then  $Q$ is {\it amenable relative to $P$} if there exists a net $\xi_n\in L^2(M)\bar{\otimes}_{P}L^2(M)$ of tracial, $Q$-almost central vectors \cite{OP07}.

\end{definition}

Thus, the failure of an algebra to be amenable (or amenable relative to some other algebra) can be viewed as a source of ``spectral gap rigidity". The notion of spectral gap rigidity has been introduced by S. Popa in the context of von Neumann algebras, and has been used to great effect in \cite{Po06a,Po06c,OP07}. 

\subsection{Deformation vs. rigidity} We finally explain how deformation and rigidity are put together in Popa's proof of his strong rigidity theorem.

\begin{theorem}[Popa, \cite{Po03,Po04}]\label{strong} Let $\Gamma$ be an icc group and $\Gamma\curvearrowright (X,\mu)=(X_0,\mu_0)^{\Gamma}$  be a Bernoulli action.  Let $\Lambda$ be a property (T) group and $\Lambda\curvearrowright (Y,\nu)$ be any free ergodic pmp action.

If $L^{\infty}(X)\rtimes\Gamma\cong L^{\infty}(Y)\rtimes\Lambda$, then the actions $\Gamma\curvearrowright (X,\mu)$ and $\Lambda\curvearrowright (Y,\nu)$ are conjugate.
\end{theorem}

In the first part of the proof, Popa essentially shows that any subalgebra $P$ of $M=L^{\infty}(X)\rtimes\Gamma$ with the  property (T) can be unitarily conjugated into $L(\Gamma)$. 
To give an idea of why this is plausible, let $\tilde M$ be the tracial  algebra containing $M$ together with its 1-parameter group of automorphisms $(\theta_t)_{t\in\mathbb R}$  constructed in Example \ref{wmall}. Then the deformation $\phi_t=E_M\circ{\theta_t}_{|M}$ converges uniformly to the identity on the unital ball of $P$, as $t\rightarrow 0$.  It is immediate that the same is true for $\theta_t$. This forces that, for small enough $t>0$, the restriction of $\theta_t$ to $P$ is inner.

If we denote $A_0=L^{\infty}(X_0)$, then $M=(\bar{\otimes}_{g\in\Gamma}A_0)\rtimes\Gamma$.
By its definition, $\theta_t$ is inner both on $L(\Gamma)$ and on $\bar{\otimes}_{g\in F}A_0$, for every finite subset $F\subset\Gamma$. As it turns out, the converse is also true: any subalgebra of $M$ on which $\theta_t$ is inner can be unitarily conjugate into one of these algebras. Since $\bar{\otimes}_{g\in F}A_0$ is abelian, it cannot contain any property (T) subalgebra. This altogether implies the initial claim.

Now, identify $L^{\infty}(X)\rtimes\Gamma=L^{\infty}(Y)\rtimes\Lambda$.  Since $L(\Lambda)$ has property (T), by the first part of the proof,  we may assume that $L(\Lambda)\subset L(\Gamma)$. In the second part of the proof, Popa proves that any group measure space Cartan subalgebra of $M$ that is normalized by many unitaries in $L(\Gamma)$  can be unitarily conjugate into $L^{\infty}(X)$. 

Thus, we have that both $L(\Lambda)\subset L(\Gamma)$ and $uL^{\infty}(Y)u^*\subset L^{\infty}(X)$, for some unitary $u\in M$. In the  final part of the proof, Popa 
is able to show that we may take $u=1$. This readily implies that the actions are conjugate.

\section{W$^*$-superrigidity for Bernoulli actions and generalized wreath product groups} 

The natural question underlying Theorem \ref{strong} is whether Bernoulli actions of icc property (T) groups are W$^*$-superrigid \cite{Po04}. In the first part of this section, we  discuss the  recent positive resolution of this question.

\begin{theorem}[Ioana, \cite{Io10}]\label{ber}
Let $\Gamma$ be any icc property (T) group.
Then the Bernoulli action $\Gamma\curvearrowright  (X,\mu)=(X_0,\mu_0)^{\Gamma}$ is W$^*$-superrigid. \end{theorem} 
More generally, Theorem \ref{ber} holds for icc groups $\Gamma$ which admit an infinite normal subgroup with the relative property (T).

To outline the strategy of proof, denote $A=L^{\infty}(X)$ and $M=A\rtimes\Gamma$. Assume that  $M=B\rtimes\Lambda$, where $B=L^{\infty}(Y)$. This new group measure space decomposition of $M$ gives rise to an embedding $\Delta:M\rightarrow M\bar{\otimes}M$ defined by $\Delta(bv_h)=bv_h\otimes v_h$, for all $b\in B$ and $h\in\Lambda$, where $(v_h)_{h\in\Lambda}\subset M$ denote the canonical unitaries \cite{PV09}.

The proof relies on a classification of all possible embeddings $\Delta:M\rightarrow M\bar{\otimes}M$. This classification is precise enough to imply that $A$ and $B$ are unitarily conjugated, and hence that the actions $\Gamma\curvearrowright (X,\mu)$ and $\Lambda\curvearrowright (Y,\nu)$ are orbit equivalent.  Since the action $\Gamma\curvearrowright (X,\mu)$ is OE superrigid by a result of Popa  \cite{Po05}  it follows that the actions are indeed conjugate.

Now, let us say a few words about the techniques that we use in order to classify embeddings $\Delta: M\rightarrow M\bar{\otimes}M$. Since the same techniques are needed in the study of embeddings
 $\theta:M\rightarrow M$, for simplicity, we only discuss the latter issue here.   Since $\Gamma$ has property (T), by the first part of the proof of Theorem \ref{strong}, we may assume that $\theta(L(\Gamma))\subset L(\Gamma)$.
 Denoting $D=\theta(A)$, we 
 have that $D$ is an abelian algebra which is normalized by a group of unitary elements, $\theta(\Gamma)$, from $L(\Gamma).$

 The main novelty of \cite{Io10} is a structural result for abelian subalgebras $D$ of $M$: suppose that $D$ is normalized by a sequence of unitary elements $u_n\in L(\Gamma)$ such that $u_n\rightarrow 0$, weakly. Then essentially  either $D$  can be unitarily conjugated into  $L(\Gamma)$, or  its relative commutant $D'\cap M$ can be unitarily conjugated into $A$. Assuming that $D$ cannot be unitarily conjugated into $L(\Gamma)$, we deduce that the unitaries $u_n\in L(\Gamma)$ are uniformly close to the discrete subgroup $\Gamma\subset L(\Gamma)$. More precisely, there exists a sequence $g_n\in\Gamma$ such that $\sup_n\|u_n-u_{g_n}\|_2<\sqrt{2}$.  This allows us to carry subsequent calculations by analogy with the case $u_n\in\Gamma$ and thereby conclude that $D'\cap M$ can be unitarily conjugate into $A$. 
 
The techniques from $\cite{Io10}$ also yield a new class of II$_1$ factors that do not arise from groups. In the setting of Theorem \ref{ber}, assume additionally that $\Gamma$ is torsion free. Then for any projection $p\in M\setminus\{0,1\}$, the II$_1$ factor $pMp$ is not  isomorphic to the group von Neumann algebra of any countable group. Moreover, $pMp$ is not isomorphic to any twisted group von Neumann algebra $L_{\alpha}(G)$, where $\alpha$ is a scalar $2$-cocycle on a countable group $G$. This gave the first examples of such II$_1$ factors.
 
 \vskip 0.1in
In the second part of this section, we present the recent discovery of the first classes of W$^*$-superrigid groups.

\begin{theorem}[Ioana, Popa, Vaes, \cite{IPV10}]\label{ipv} Let $\Gamma_0$ be any non-amenable group and $S$ be any infinite amenable group. Let $\Gamma_0^{(S)}=\oplus_{s\in S}\Gamma_0$ and define the wreath product group $\Gamma=\Gamma_0^{(S)}\rtimes S$.  Consider the left multiplication action of $\Gamma$ on the coset space $I=\Gamma/S$.

Then the generalized wreath product $G={({\mathbb Z}/{2\mathbb Z})}^{(I)}\rtimes\Gamma$ is W$^*$-superrigid. That is, if $\Lambda$ is any countable group such that $L(G)\cong L(\Lambda)$, then $G\cong\Lambda$.

\end{theorem}  
 
The fact that certain generalized wreath product groups are W$^*$-superrigid should not be surprising, since such groups have been recently recognized to be remarkably rigid in the von Neumann algebra context. For instance, Popa's strong rigidity theorem implies that  if $G_i={({\mathbb Z}/{2\mathbb Z})}^{(\Gamma_i)}\rtimes\Gamma_i$, where $\Gamma_i$ is a property (T) group for $i\in\{1,2\}$, then $L(G_1)\cong L(G_2)$ entails $G_1\cong G_2$ \cite{Po04}. 

Note, however, that the conclusion of Theorem \ref{ipv} does not hold for plain wreath product groups $G={({\mathbb Z}/{2\mathbb Z})}^{(\Gamma)}\rtimes\Gamma$. In fact, for any non-trivial torsion free group $\Gamma$, there exists a torsion free group $\Lambda$ such that $L(G)\cong L(\Lambda)$ \cite{IPV10}. Nevertheless, for certain classes of groups $\Gamma$, including icc property (T) groups and products of non-amenable groups, we are still able to classify more or less explicitly all groups $\Lambda$ with $L(G)\cong L(\Lambda)$.

To describe the main steps of the proof of Theorem \ref{ipv}, let $M=L(G)$ and assume that $M=L(\Lambda)$, for a countable group $\Lambda$. Denote by $\Delta_{\Lambda}:M\rightarrow M\bar{\otimes}M$ the embedding given by $\Delta_{\Lambda}(v_h)=v_h\otimes v_h$, where $(v_h)_{h\in\Lambda}\subset M$ are the canonical unitaries.   Similarly, define $\Delta_G:M\rightarrow M\bar{\otimes}M$. We start by viewing $M$ as the group measure space II$_1$ factor of the generalized Bernoulli action $\Gamma\curvearrowright (\{0,1\},\mu_0)^{\Gamma}$, where $\mu_0$ is the measure on $\{0,1\}$ given by $\mu_0(\{0\})=\mu_0(\{1\})=\frac{1}{2}$. 
By extending the methods of \cite{Io10} from plain Bernoulli actions to generalized Bernoulli actions, we then give a classification of all  possible embeddings $\Delta:M\rightarrow M\bar{\otimes}M$. 

When applied to $\Delta_{\Lambda}$, this enables us to deduce the existence of a unitary element $\Omega\in M\bar{\otimes}M$ such that $\Delta_{\Lambda}(x)=\Omega\Delta_{G}(x)\Omega^*$, for all $x\in M$.  Moreover, it follows that $\Omega$ satisfies a certain ``dual" 2-cocycle relation. A main novelty of \cite{IPV10} is a vanishing result for dual 2-cocycles which allows to conclude that the groups $G$ and $\Lambda$ are isomorphic from the existence of $\Omega$.
  
Let us emphasize a particular case of this result that provides a surprising criterion for the unitary conjugacy of arbitrary icc groups $G,\Lambda$ giving the same II$_1$ factor, $L(G)=L(\Lambda)$. Assume that there exists a constant $\kappa<\sqrt{2}$ with the property that for every $g\in G$ we can find $h\in\Lambda$ such that $\|u_g-v_h\|_2\leqslant\kappa$. Then $G\cong\Lambda$ and there exist a group isomorphism $\delta:G\rightarrow\Lambda$, a character $\eta:G\rightarrow\mathbb T$ and a unitary element $u\in L(G)$ such that $uu_gu^*=\eta(g)v_{\delta(g)}$, for all $g\in G$.

\section{Uniqueness of group measure space Cartan subalgebras}
In this section we discuss several uniqueness results for group measure space Cartan subalgebras. We start with a  general uniqueness result of Popa and Vaes:

\begin{theorem}[Popa, Vaes, \cite{PV09}]\label{pv} Let $\Gamma=\Gamma_1*_{\Sigma}\Gamma_2$ be a non-trivial amalgamated free product such that $\Gamma_1$ admits a non-amenable subgroup  with the relative property (T), $\Sigma$ is amenable and  there exist $g_1,g_2,...,g_n\in\Gamma$ such that $\cap_{i=1}^ng_i\Sigma g_i^{-1}$ is finite.
Let $\Gamma\curvearrowright (X,\mu)$ be any free ergodic pmp action.

Then  $L^{\infty}(X)\rtimes\Gamma$ has a unique group measure space Cartan subalgebra, up to unitary conjugacy.
\end{theorem}

Theorem \ref{pv} covers in particular arbitrary free ergodic pmp actions of any free product  $\Gamma=\Gamma_1*\Gamma_2$ of an infinite property (T) group and a non-trivial group.
In combination with the OE superrigidity theorems of Popa \cite{Po05,Po06a} and Kida \cite{Ki09}, it lead in \cite{PV09} to the first families of W$^*$-superrigid actions. For instance,  if $T_n<\text{PSL}_n(\mathbb Z)$ is the group of triangular matrices for some $n\geqslant 3$, then any free {\it mixing} pmp action of $\text{PSL}_n(\mathbb Z)*_{T_n}\text{PSL}_n(\mathbb Z)$ is W$^*$-superrigid. For a group whose all free ergodic pmp actions are W$^*$-superrigid, see \cite{HPV10}.

To give an overview of the proof, assume for simplicity that $\Sigma=\{e\}$.
Define $A=L^{\infty}(X)$ and $M=A\rtimes\Gamma$. Then we have an amalgamated free product decomposition $M=M_1*_{A}M_2$, where $M_1=A\rtimes\Gamma_1$ and $M_2=A\rtimes\Gamma_2$. Let $(\tilde M, (\theta_t)_{t\in\mathbb R})$ be the malleable deformation constructed in Example \ref{ipp}.

Assume that $M=B\rtimes\Lambda$ is another group measure space decomposition and denote by $(v_h)_{h\in\Lambda}$ the canonical unitaries.  The first part of the proof amounts to transferring some of the rigidity of $\Gamma$ to $\Lambda$. 
Intuitively, since $\Gamma$ has a non-amenable subgroup with the relative property (T) while $A$ and $B$ are amenable, $\Lambda$ must admit a ``non-amenable subset $S$ with the relative property (T)". 
Concretely, Popa and Vaes prove that given $\varepsilon>0$, there exist a sequence $h_n\in\Lambda$ and $t>0$ such that the unitary elements $v_n:=v_{h_n}$ satisfy \begin{enumerate}
\item\label{rigid} $\|\theta_t(v_{n})-v_{n}\|_2<\varepsilon$, for all $n$, and 
\item\label{mixing} $\|E_A(xv_{n}y)\|_2\rightarrow 0$, for all $x,y\in M$.
\end{enumerate}

Condition  (1) implies that the unitaries $v_{n}$ have ``uniformly bounded length" in $M={M_1}*_{A}M_2$: they are almost supported on words on length $\leqslant\kappa$ in $M_1$ and $M_2$, for a fixed $\kappa\geqslant 1$. In the second part of the proof, using a combinatorial argument, Popa and Vaes  conclude that since $B$ is abelian and is normalized by the unitaries $v_{n}$, the whole unit ball of $B$ must also have  uniformly bounded length. It is clear that both $M_1$ and $M_2$ have uniformly bounded length. Conversely, the main technical result of \cite{IPP05} shows that any subalgebra $B$ with this property can be unitarily conjugated into either $M_1$ or $M_2$. Finally, since the normalizer of $B$ generates the whole $M$, this forces that $B$ can be unitarily conjugate into $A$.

\vskip 0.1in
Next, we comment on a result of Chifan and Peterson which gives a 1-cohomology approach to uniqueness of group measure space Cartan subalgebras.

\begin{theorem}[Chifan, Peterson, \cite{CP10}]\label{cp} Let $\Gamma$ be a countable group which admits a non-amenable subgroup with the relative property (T) and an unbounded cocycle $c:\Gamma\rightarrow H$ into a mixing orthogonal representation $\pi:\Gamma\rightarrow\mathcal O(H)$. Let $\Gamma\curvearrowright (X,\mu)$ be any free ergodic pmp action.

Then $L^{\infty}(X)\rtimes\Gamma$ has a unique group measure space Cartan subalgebra, up to unitary conjugacy.
\end{theorem}

Recall that  $\pi$ is mixing if $\langle\pi(g)\xi,\eta\rangle\rightarrow 0$, as $g\rightarrow\infty$, for any vectors $\xi,\eta\in H$. 

Define $A=L^{\infty}(X)$ and $M=A\rtimes\Gamma$. The original proof of Theorem \ref{cp} uses Peterson's technique of closable derivations  \cite{Pe06}. Following \cite{Va10b}, we consider instead the deformation of $M$ given by the cocycle $c:\Gamma\rightarrow H$ as in Example \ref{si}.

The proof of Theorem \ref{cp} has the same skeleton as the proof of Theorem \ref{pv}. For this reason, we only emphasize two new ingredients.  Thus, let $B\subset M$ be an abelian von Neumann subalgebra. Assume for every $\varepsilon>0$ there exist $t>0$ and a sequence of unitaries $v_n$ normalizing $B$ such that conditions (1) and (2) hold. Chifan and Peterson then prove that the deformation $\theta_t$ must converge uniformly on the unit ball of $B$. If this is the case, then a result from \cite{Pe06} further implies that either $B$ can be unitarily conjugate into $A$, or $\theta_t$ converges uniformly on the normalizer of $B$. When $B$ is Cartan subalgebra, the normalizer of $B$ generates $M$, and the latter condition is impossible, since $c$ is unbounded.

 \vskip 0.05in 
 
In \cite{Va10b},  Vaes generalized Theorem \ref{cp} by replacing the condition  that $\pi$ is mixing with the weaker condition that $\pi$ is mixing relative to a family of amenable subgroups of $\Gamma$. This result also recovers Theorem \ref{pv} since any amalgamated free product  $\Gamma=\Gamma_1*_{\Sigma}\Gamma_2$ admits an unbounded cocycle into the quasi-regular representation $\pi:\Gamma\rightarrow\ell^2(\Gamma/\Sigma)$ which is mixing relative to $\Sigma$. 
 \vskip 0.05in
Theorems \ref{pv} and \ref{cp}  provide supportive evidence for the general conjecture that $L^{\infty}(X)\rtimes\Gamma$ must have a unique Cartan subalgebra, for any free ergodic pmp action of any group with $\beta_1^{(2)}(\Gamma)>0$. 
 The following result provides further positive evidence towards this conjecture.

\begin{theorem} [Ioana, \cite{Io11a,Io11b}]\label{io} Let $\Gamma$ be a countable group with $\beta_1^{(2)}(\Gamma)>0$. 
Let $\Gamma\curvearrowright (X,\mu)$ be a free ergodic pmp action which is either rigid or profinite.

Then $L^{\infty}(X)\rtimes\Gamma$ has a unique group measure space Cartan subalgebra, up to unitary conjugacy.
\end{theorem}

Recall that an action $\Gamma\curvearrowright (X,\mu)$ is {\it rigid} if the inclusion $L^{\infty}(X)\subset L^{\infty}(X)\rtimes\Gamma$ has the relative property (T) \cite{Po01b}. Examples of rigid actions are given by the actions $\text{SL}_2(\mathbb Z)\curvearrowright\mathbb T^2$ \cite{Po01b} and SL$_2(\mathbb Z)\curvearrowright \text{SL}_2(\mathbb R)/\text{SL}_2(\mathbb Z)$ \cite{IS10}. By \cite{Ga08} any free product group $\Gamma=\Gamma_1*\Gamma_2$ with $|\Gamma_1|\geqslant 2$ and $|\Gamma_2|\geqslant 3$ admits a continuum of rigid actions whose II$_1$ factors are mutually non-isomorphic.
 
 Also, recall that an action $\Gamma\curvearrowright (X,\mu)$ is {\it profinite} if it is the inverse limit $\varprojlim\Gamma\curvearrowright (X_n,\mu_n)$ of  actions of $\Gamma$ on finite probability spaces $(X_n,\mu_n)$. Note that if  $G=\varprojlim\Gamma/\Gamma_n$ is the profinite completion of $\Gamma$ with respect to a descending chain of  finite index subgroups, then the left translation action $\Gamma\curvearrowright G$ is profinite.
 
 To outline  the proof of Theorem \ref{io}, define $A=L^{\infty}(X)$ and $M=A\rtimes\Gamma$.  Suppose that $M=B\rtimes\Lambda$ is another group measure space decomposition.
 
In the first part of the proof, we show  $A$ can be unitarily conjugated into $B\rtimes\Sigma$, for some amenable subgroup $\Sigma<\Lambda$.  This is quite unexpected, because  a priori one has no knowledge about the subgroups of  $\Lambda$. To derive this, we consider the deformation of $M$ arising from an unbounded cocycle $c:\Gamma\rightarrow\ell^2(\Gamma)$ (such a cocycle exists since $\beta_1^{(2)}(\Gamma)>0$ \cite{PT07}). We then combine the above mentioned results of \cite{CP10} with quite delicate estimates in an ultraproduct algebra $M^{\mathcal \omega}$ associated with a cofinal ultrafilter $\omega$ over a (possibly uncountable) directed set. 
  
 Since $\Sigma$ is amenable, the algebra $B\rtimes\Sigma$ is also amenable. Thus, we may essentially assume that von Neumann algebra $N$ generated by $A$ and $B$ is amenable. In the second part of the proof,  we use this to conclude that $A$ and $B$ are unitarily conjugate.
If $A$ and $B$ are not conjugate, then the equivalence relation $R$ on $(X,\mu)$ associated with the inclusion $A\subset N$ \cite{FM75} must be ``weakly normal" in $R_{\Gamma\curvearrowright X}$. On the other hand, since $\Gamma$ has positive first $\ell^2$-Betti number,  $R_{\Gamma\curvearrowright X}$ does not admit an aperiodic,  amenable subequivalence relation that is weakly normal  \cite{Ga99,Ga01}. 
 
\section{Uniqueness of arbitrary Cartan subalgebras}
In this section we discuss several results providing classes of II$_1$ factors with a unique Cartan subalgebra, up to unitary conjugacy. The proofs of these results make crucial use of the following approximation property for groups introduced by Cowling and Haagerup:

\begin{definition}
A countable group $\Gamma$ is  {\it weakly amenable} \cite{CH88} if there exists a sequence a sequence of functions  $\varphi_k:\Gamma\rightarrow\mathbb C$   such that 
\begin{itemize}
\item $\varphi_k$ has finite support, for all $k$,

\item $\lim_k\varphi_k(g)=1$, for all $g\in\Gamma$,

\item $\limsup_{k}\|\phi_k\|_{\text{cb}}<\infty$, where $\phi_k:L(\Gamma)\rightarrow L(\Gamma)$ is the  unique map satisfying $\phi_k(u_g)=\varphi_k(g)u_g$, for all $g\in\Gamma$, and $\|\phi_k\|_{\text{cb}}$ is its completely bounded norm.
\end{itemize}
Moreover, if  there exist  $\varphi_k:\Gamma\rightarrow\mathbb C$ as above such that $\limsup_{k}\|\phi_k\|_{\text{cb}}=1$, then we say that $\Gamma$ has the {\it complete metric approximation property} (CMAP) \cite{Ha79}.

\end{definition}

The first result ever showing uniqueness, up to unitary conjugacy, of arbitrary Cartan subalgebras,  was obtained by Ozawa and Popa:

\begin{theorem}
[Ozawa, Popa, \cite{OP07}]\label{op} Let $\mathbb F_n\curvearrowright (X,\mu)$ be a free ergodic profinite pmp action of a free group, $\mathbb F_n$, for some $n\geqslant 2$.

Then $L^{\infty}(X)\rtimes\Gamma$ has a unique Cartan subalgebra, up to unitary conjugacy.
\end{theorem}

Denote $A=L^{\infty}(X)$ and $M=A\rtimes\mathbb F_n$.  Since $\mathbb F_n$ has the CMAP \cite{Ha79} and the action $\mathbb F_n\curvearrowright X$ is profinite,  the II$_1$ factor $M$ also has the CMAP: there exists a sequence of finite rank completely bounded maps $\phi_k:M\rightarrow M$ such that $\|\phi_k(x)-x\|_2\rightarrow 0$, for all $x\in M$, and $\limsup_{k}\|\phi_k\|_{\text{cb}}=1$. 
Consider an arbitrary diffuse amenable von Neumann subalgebra $P\subset M$ and denote by $\mathcal G$ its normalizer. 

Ozawa and Popa made the amazing discovery that since $M$ has the CMAP, 
the action of $\mathcal G$ on $P$ by conjugation is {\it weakly compact}.
More precisely, there exists a net of positive vectors $\xi_k\in L^2(P)\bar{\otimes}L^2(P)$ which are tracial, $P$-almost central and almost invariant under the diagonal action of $\mathcal G$. Note that if the vectors $\xi_k$ are actually invariant under the diagonal action of $\mathcal G$, then  the action $\mathcal G\curvearrowright P$ is {\it compact}, i.e. the closure of $\mathcal G$ inside Aut$(P)$ is compact.

In the second part of the proof, Ozawa and Popa combine the free malleable deformation of $M$ \cite{Po86,Po06c} with the weak compactness of the action $\mathcal G\curvearrowright P$. They conclude that either a corner of $P$ embeds into $A$ or the von Neumann algebra generated by $\mathcal G$ is amenable. If $P\subset M$ is a Cartan subalgebra, then $\mathcal G''=M$ is not amenable. Therefore, by Theorem \ref{inter}, $P$ must be unitarily conjugate to $A$.

\vskip 0.05in

In \cite{Oz10} Ozawa showed that one can replace the usage of the CMAP by that of weak amenability in  the proof of Theorem \ref{op}. This result opened up the possibility that Theorem \ref{op} could be extended to weakly amenable groups $\Gamma$ that do not have the CMAP. 
Motivated by this, Chifan and Sinclair proved the following:

\begin{theorem}
[Chifan, Sinclair, \cite{CS11}]\label{cs} Let $\Gamma\curvearrowright (X,\mu)$ be a free ergodic profinite pmp action of a non-elementary hyperbolic group $\Gamma$.

Then $L^{\infty}(X)\rtimes\Gamma$ has a unique Cartan subalgebra, up to unitary conjugacy.
\end{theorem}

 Without explaining further details, let us mention the three main ingredients of the proof of Theorem \ref{cs}. Firstly, since hyperbolic groups are weakly amenable \cite{Oz07}, the conjugation action $\mathcal N_{M}(P)\curvearrowright P$ is weakly compact, for any diffuse von Neumann subalgebra $P$ of  $M=L^{\infty}(X)\rtimes\Gamma$ \cite{Oz10}.
Secondly, following \cite{MMS04}, for any hyperbolic  group $\Gamma$ there is a proper quasi-cocycle $c:\Gamma\rightarrow\ell^2(\Gamma)$.  
This gives rise to a  C$^*$-algebraic deformation of $M$ (\cite{CS11}, see Example \ref{chsi}) that is compact relative to $L^{\infty}(X)$.  Finally, since hyperbolic groups are exact,  elements of $M$ admit ``good"  approximations by elements in the  reduced C$^*$-algebra  $L^{\infty}(X)\rtimes_{\text{r}}\Gamma$.

\vskip 0.05in
Very recently, Popa and Vaes obtained sweeping ``unique Cartan subalgebra" results. In particular, they were able to extend Theorems \ref{op} and \ref{cs} to {\it arbitrary} actions of free groups and non-elementary hyperbolic groups. 

\begin{theorem}[Popa, Vaes, \cite{PV11,PV12}]\label{pv11} Let $\Gamma$ be a countable group and  let $\Gamma\curvearrowright (X,\mu)$ be any free ergodic pmp action. Assume that either
\begin{enumerate}
\item
 $\Gamma$ is weakly amenable and admits an unbounded cocycle into a  non-amenable mixing orthogonal representation $\pi:\Gamma\rightarrow\mathcal O(H)$, or  
 \item$\Gamma$ is non-elementary hyperbolic. 
 
 Then $L^{\infty}(X)\rtimes\Gamma$ has a unique Cartan subalgebra, up to unitary conjugcy.
\end{enumerate}
\end{theorem} 

Recall that a representation $\pi$ is amenable if $\pi\otimes\bar{\pi}$ has almost invariant vectors. Hence, the left regular representation of a non-amenable group $\Gamma$ is non-amenable. Therefore, any weakly amenable group $\Gamma$ with $\beta_1^{(2)}(\Gamma)>0$ satisfies Theorem \ref{pv11}.

Theorem \ref{pv11} has the following beautiful consequence.
If $2\leqslant m,n\leqslant \infty$ and $m\not=n$, then any free ergodic pmp actions  $\mathbb F_m\curvearrowright (X,\mu)$ and $\mathbb F_n\curvearrowright (Y,\nu)$ give rise to non-isomorphic II$_1$ factors, $L^{\infty}(X)\rtimes\mathbb F_m\not\cong L^{\infty}(Y)\rtimes\mathbb F_n.$

 Indeed, if these factors were isomorphic, then by Theorem \ref{pv11}, the actions $\mathbb F_m\curvearrowright X$ and $\mathbb F_n\curvearrowright Y$ would  be orbit equivalent. However, it is proven in \cite{Ga99,Ga01} that 
free actions of free groups of different ranks are never orbit equivalent.

In combination with \cite{Bo09a,Bo09b}, Theorem \ref{pv11} also allows to completely classify all amplifications of the II$_1$ factors associated with the wreath product groups $\mathbb Z\wr\mathbb F_n$. Thus, $L(\mathbb Z\wr\mathbb F_m)^t\cong L(\mathbb Z\wr\mathbb F_n)^s$ if and only if $(m-1)/s=(n-1)/t$.

In order to give an overview the proof of Theorem \ref{pv11}, denote $M=L^{\infty}(X)\rtimes\Gamma$.  Since $\Gamma$ is weakly amenable, there exists a sequence $\phi_k:M\rightarrow M$ of  completely bounded maps  of ``finite rank relative to $L^{\infty}(X)$" such that $\|\phi_k(x)-x\|_2\rightarrow 0$, for all $x\in M$, and  $\limsup_{k}\|\phi_k\|_{\text{cb}}=1$. The existence of such maps, however, cannot imply that actions of the form $\mathcal N_M(P)\curvearrowright P$ are weakly compact.  Indeed, the action $\Gamma\curvearrowright X$ might itself not be weakly compact (e.g. if it is a Bernoulli action and $\Gamma$ is non-amenable). 

Nevertheless, Popa and Vaes discovered that there is an appropriate notion of weak compactness  which is some sense relative to $L^{\infty}(X)$. They then proved, by extending a technique from \cite{OP07}, that  the action $\mathcal N_M(P)\curvearrowright P$ has the relative weak compactness property, for any diffuse amenable subalgebra $P\subset M$.
To complete the proof, Popa and Vaes combine this  property with the malleable deformation of $M$ associated with an unbounded cocycle $c:\Gamma\rightarrow H$ \cite{Si10},  in case (1), and a proper quasi-cocycle $c:\Gamma\rightarrow\ell^2(\Gamma)$ \cite{CS11}, in case (2). As shown in a later version of \cite{PV12}, in case (2) one can alternatively use the fact that hyperbolic groups are biexact  \cite{Oz03}.

\section{References}

\vskip 0.2in
\renewcommand{\refname}{}    
\vspace*{-36pt}              

\end{document}